\documentclass[11pt]{article}

\usepackage[a4paper,margin=1.08in]{geometry}
\usepackage[T1]{fontenc}
\usepackage[utf8]{inputenc}
\usepackage{lmodern}
\usepackage{microtype}
\usepackage{amsmath,amssymb,amsthm,mathtools}
\usepackage{enumitem}
\usepackage{hyperref}
\usepackage{aliascnt}
\usepackage[nameinlink,noabbrev]{cleveref}
\usepackage{orcidlink}

\hypersetup{
  hidelinks,
  pdftitle={Fractional minimality of the Simons cone for s close to 1},
  pdfauthor={Matteo Novaga and Emanuele Paolini}
}
\numberwithin{equation}{section}
\setlist[enumerate]{leftmargin=2.2em,itemsep=0.25em,topsep=0.4em}

\newtheorem{theorem}{Theorem}[section]

\newaliascnt{proposition}{theorem}
\newtheorem{proposition}[proposition]{Proposition}
\aliascntresetthe{proposition}

\newaliascnt{lemma}{theorem}
\newtheorem{lemma}[lemma]{Lemma}
\aliascntresetthe{lemma}

\newaliascnt{corollary}{theorem}
\newtheorem{corollary}[corollary]{Corollary}
\aliascntresetthe{corollary}

\newaliascnt{definition}{theorem}
\newtheorem{definition}[definition]{Definition}
\aliascntresetthe{definition}

\newaliascnt{remark}{theorem}
\newtheorem{remark}[remark]{Remark}
\aliascntresetthe{remark}

\crefname{theorem}{Theorem}{Theorems}
\crefname{proposition}{Proposition}{Propositions}
\crefname{lemma}{Lemma}{Lemmas}
\crefname{corollary}{Corollary}{Corollaries}
\crefname{definition}{Definition}{Definitions}
\crefname{remark}{Remark}{Remarks}

\DeclareMathOperator{\PV}{PV}
\newcommand{\R}{\mathbb R}
\newcommand{\Sph}{\mathbb S}
\newcommand{\1}{\mathbf 1}
\newcommand{\dd}{\,d}
\newcommand{\cQ}{\mathcal Q}
\newcommand{\cH}{\mathcal H}
\newcommand{\cD}{\mathcal D}
\DeclareMathOperator{\sgn}{sgn}
\newcommand{\eps}{\varepsilon}

\title{\textbf{Fractional minimality of the Simons cone \\ for $s$ close to $1$}}

\author{
  Matteo Novaga\thanks{Dipartimento di Matematica, Universit\`a di Pisa,
  Largo Bruno Pontecorvo 5, 56127 Pisa, Italy.
  Email: \href{mailto:matteo.novaga@unipi.it}{matteo.novaga@unipi.it}.
  ORCID: \orcidlink{0000-0002-0812-3780}}
  \and
  Emanuele Paolini\thanks{Dipartimento di Matematica, Universit\`a di Pisa,
  Largo Bruno Pontecorvo 5, 56127 Pisa, Italy.
  Email: \href{mailto:emanuele.paolini@unipi.it}{emanuele.paolini@unipi.it}.
  ORCID: \orcidlink{0000-0001-7021-7537}}
}

\date{}

\begin{document}
\maketitle

\begin{abstract}
\noindent Let $\Sigma:=\{(X,Y)\in\R^4\times\R^4:|X|=|Y|\}$
be the Simons cone and let $C:=\{|X|\le |Y|\}$. 
We prove that $C$ minimizes the fractional perimeter $P_s$, under
compact perturbations, whenever $s$ is sufficiently close to $1$.  The proof
constructs two homogeneous families of strict one-sided barriers and converts
them into nonlocal sub-calibrations.  
\end{abstract}



\section{Introduction}\label{sec:introduction}

The Simons cone is the singular hypersurface
\begin{equation}\label{eq:Sigma-intro}
 \Sigma=\{(X,Y)\in\R^4\times\R^4:|X|=|Y|\}\subset\R^8.
\end{equation}
We denote by
\begin{equation}\label{eq:C-intro}
 C:=\{(X,Y)\in\R^4\times\R^4:|X|\le |Y|\}
\end{equation}
one of the two closed cones bounded by $\Sigma$.  Simons proved the relevant
stability and low-dimensional rigidity results for minimal cones
\cite{Simons}; Bombieri, De Giorgi and Giusti subsequently proved that the cone is area minimizing, equivalently that $C$ is locally perimeter minimizing, \cite{BDG}.  Calibration and
sub-calibration proofs were later given, among others, by Davini
\cite{Davini} and by De Philippis and Paolini \cite{DePhilippisPaolini}.

We study the analogous minimizing property for the fractional perimeter
introduced by Caffarelli, Roquejoffre and Savin \cite{CRS}.  For an open set
$\Omega\subset\R^n$, $s\in(0,1)$ and a measurable set $E\subset\R^n$, we use
the normalization
\begin{align}\label{eq:perimeter-intro}
 P_s(E,\Omega)
 &:={\frac12}\int_{\Omega}\int_{\Omega}
 \frac{\lvert\chi_E(x)-\chi_E(y)\rvert}{|x-y|^{n+s}}\dd x\dd y \\
 &\quad+\int_{\Omega}\int_{\R^n\setminus\Omega}
 \frac{\lvert\chi_E(x)-\chi_E(y)\rvert}{|x-y|^{n+s}}\dd x\dd y.\nonumber
\end{align}
At a smooth boundary point $p\in\partial E$, its Euler--Lagrange operator is
the fractional mean curvature
\begin{equation}\label{eq:Hs-intro}
 H_s^E(p):=\PV\int_{\R^n}
 \frac{\chi_{E^c}(q)-\chi_E(q)}{|q-p|^{n+s}}\dd q.
\end{equation}
The fractional perimeter converges, after the usual normalization, to the
classical perimeter as $s\uparrow1$ \cite{CaffarelliValdinoci}; the corresponding
curvatures converge to the classical mean curvature on smooth families
\cite{AbatangeloValdinoci,CesaroniDeLucaNovagaPonsiglione}.

Nonlocal minimal cones play a central role in the regularity theory.  In the
plane, halfplanes are the only minimizing cones \cite{SavinValdinoci2D}.
D\'avila, del Pino and Wei constructed nonlocal Lawson cones and studied their
stability \cite{DavilaDelPinoWei}.  Recently, Palatucci \cite{Palatucci}
exhibited a singular Lawson cone minimizing the fractional perimeter in
dimension seven, a phenomenon with no classical counterpart, since by
Simons' theorem there are no singular area-minimizing cones in dimension at
most seven \cite{Simons}.  In the symmetric $4+4$ case, the exchange
map $(X,Y)\mapsto(Y,X)$ shows directly that the Simons cone is stationary for
every $s\in(0,1)$.  Stationarity, however, does not imply minimality; the
purpose of this paper is to establish the latter in the perturbative regime
$s\uparrow1$.

The proof is based on nonlocal calibrations.  Cabr\'e developed calibrations
and null-Lagrangians associated with foliations by solutions of the nonlocal
mean-curvature equation \cite{Cabre}.  Pagliari introduced an antisymmetric
scalar calibration for nonlocal total-variation energies and used it to prove
the minimality and uniqueness of halfspaces \cite{Pagliari}.  Here the leaves
are strict sub- and supersolutions rather than exact solutions, so the
appropriate object is a one-sided nonlocal sub-calibration.

The main difficulty is uniformity at infinity.  Set
\[
 E_1^+:=\{|X|^4\le |Y|^4+1\}.
\]
On compact portions of $\partial E_1^+$, positivity of $H_s^{E_1^+}$ follows
from the local limit $s\uparrow1$ and positivity of the classical mean
curvature.  This argument alone is not uniform at infinity, because the
classical mean curvature decays like $y^{-13}$ along the boundary.  Rescaling
by $\eta=y^{-4}$ transforms the remote boundary point into
\[
 P_\eta=((1+\eta)^{1/4}e_1,e_1)
\]
and the set into
\[
 E_\eta=\{|X|^4\le |Y|^4+\eta\}.
\]
The key estimate is
\begin{equation}\label{eq:key-estimate-intro}
 H_s^{E_\eta}(P_\eta)\ge A_s\eta-C\eta^2,
 \qquad A_s\ge c>0,
\end{equation}
uniformly for $s$ in a one-sided neighborhood of $1$.  A Laplace transform
turns the curvature into a weighted integral of Gaussian sublevel
probabilities.  The first variation is computed explicitly using modified
Bessel functions.  The remainder requires two regimes: an exact noncentral
chi-square representation for bounded Gaussian scale, and a quantitative
normal-graph expansion with two $\eta$-derivatives for large Gaussian scale.
The latter isolates a nonnegative local mean-curvature term and leaves an
integrable $O(\eta^2t^{-3/2})$ error.

The estimate \eqref{eq:key-estimate-intro} yields strict mean-convexity of all
positive homogeneous leaves of $\Phi(X,Y)=|X|^4-|Y|^4$ and, by symmetry, of
the corresponding leaves of $-\Phi$.  The two resulting one-sided
calibrations give outward and inward minimality of $C$; submodularity then
produces minimality against arbitrary compact perturbations.

We now state the result.  Throughout the paper,
\begin{equation}\label{eq:Phi}
 \Phi(X,Y):=|X|^4-|Y|^4,
 \qquad C:=\{\Phi\le0\}.
\end{equation}
Sets are identified up to Lebesgue-null modifications.

\begin{theorem}\label{thm:main}
There exists $s_0\in(0,1)$ such that, for every $s\in(s_0,1)$, the set $C$ is
minimizing for the fractional perimeter under compact perturbations.  More
precisely, for every bounded open set $\Omega\subset\R^8$ and every measurable
set $F\subset\R^8$ satisfying $F\mathbin\Delta C\Subset\Omega$,
\[
 P_s(C,\Omega)\le P_s(F,\Omega).
\]
\end{theorem}

The analytical core is the following barrier theorem.

\begin{theorem}\label{thm:barrier}
There exists $s_0\in(0,1)$ such that, for every $s\in(s_0,1)$,
\begin{equation}\label{eq:barrier-positive}
 H_s^{E_1^+}(p)>0\qquad\text{for every }p\in\partial E_1^+,
\end{equation}
where
\[
 E_1^+:=\{|X|^4\le |Y|^4+1\}.
\]
By the symmetry $(X,Y)\mapsto(Y,X)$, the analogous statement holds for
\[
 E_1^-:=\{|Y|^4\le |X|^4+1\}.
\]
\end{theorem}

\begin{remark}\label{rem:sign}
We use
\begin{equation}\label{eq:Hs-def}
 H_s^E(p):=\PV\int_{\R^8}
 \frac{\chi_{E^c}(q)-\chi_E(q)}{|q-p|^{8+s}}\dd q.
\end{equation}
With this convention Euclidean balls have positive fractional mean curvature.
\end{remark}

\begin{remark}\label{rem:stationarity}
For every $s\in(0,1)$,
\[
 H_s^C(p)=0\qquad\text{for }p\in\Sigma\setminus\{0\}.
\]
Indeed, after applying an element of $O(4)\times O(4)$, one may assume
$p=(re_1,re_1)$.  The exchange isometry $(X,Y)\mapsto(Y,X)$ fixes $p$, swaps
$C$ and $C^c$, and preserves the kernel in \eqref{eq:Hs-def}; the principal
value therefore cancels.  Theorem~\ref{thm:main} is thus a genuine minimizing
result, not merely an Euler--Lagrange statement.
\end{remark}

The paper is organized as follows.  Section~\ref{sec:calibration} establishes
the one-sided calibration principle.  Section~\ref{sec:curvature} proves the
strict barrier theorem, including the Gaussian expansion at infinity.
Section~\ref{sec:final-summary} combines the two one-sided calibrations and
submodularity to prove Theorem~\ref{thm:main}.

\section{Calibrations and preliminary results}\label{sec:calibration}

For measurable sets $A,B\subset\R^n$, define
\[
 L_s(A,B):=\int_A\int_B\frac{1}{|x-y|^{n+s}}\dd y\dd x.
\]
For an open set $\Omega\subset\R^n$, the localized fractional perimeter is
\begin{align}\label{eq:localized-perimeter}
 P_s(E,\Omega)
 &:={L_s(E\cap\Omega,E^c\cap\Omega)}
   +{L_s(E\cap\Omega,E^c\setminus\Omega)}
   +{L_s(E\setminus\Omega,E^c\cap\Omega)}.
\end{align}
Equivalently, with
\[
 \cQ_\Omega:=\R^{2n}\setminus(\Omega^c\times\Omega^c),
\]
one has
\begin{equation}\label{eq:perimeter-Q}
 P_s(E,\Omega)=\frac12\int_{\cQ_\Omega}
 \frac{|\chi_E(x)-\chi_E(y)|}{|x-y|^{n+s}}\dd x\dd y.
\end{equation}

\begin{lemma}\label{lem:scaling}
For every $\lambda>0$,
\[
 P_s(\lambda E,\lambda\Omega)=\lambda^{n-s}P_s(E,\Omega),
 \qquad
 H_s^{\lambda E}(\lambda p)=\lambda^{-s}H_s^E(p).
\]
\end{lemma}

\begin{proof}
Both identities follow by the change of variables $x=\lambda x'$ and
$y=\lambda y'$.
\end{proof}

\begin{lemma}\label{lem:lattice}
For all measurable $A,B\subset\R^n$ and every open $\Omega$,
\[
 P_s(A^c,\Omega)=P_s(A,\Omega),
\]
and
\begin{equation}\label{eq:submod}
 P_s(A\cup B,\Omega)+P_s(A\cap B,\Omega)
 \le P_s(A,\Omega)+P_s(B,\Omega).
\end{equation}
\end{lemma}

\begin{proof}
The complement identity follows immediately from \eqref{eq:perimeter-Q}.  For
\eqref{eq:submod}, apply the pointwise lattice inequality
\[
 |(a\vee b)-(a'\vee b')|+|(a\wedge b)-(a'\wedge b')|
 \le |a-a'|+|b-b'|
\]
with $a=\chi_A(x)$, $b=\chi_B(x)$, $a'=\chi_A(y)$, $b'=\chi_B(y)$, multiply by the
kernel and integrate over $\cQ_\Omega$.
\end{proof}

\begin{definition}\label{def:one-sided-minimality}
Let $E\subset\R^n$ be measurable and let $0<s<1$.  We say that $E$ is
\emph{outward minimizing for $P_s$} if, for every bounded open set
$\Omega\subset\R^n$ and every measurable competitor $F\subset\R^n$ such that
\[
 E\subset F,
 \qquad
 F\setminus E\Subset\Omega,
\]
one has
\[
 P_s(E,\Omega)\le P_s(F,\Omega).
\]
We say that $E$ is \emph{inward minimizing for $P_s$} if, for every bounded open
set $\Omega\subset\R^n$ and every measurable competitor $F\subset\R^n$ such that
\[
 F\subset E,
 \qquad
 E\setminus F\Subset\Omega,
\]
one has
\[
 P_s(E,\Omega)\le P_s(F,\Omega).
\]
\end{definition}

We shall use the following one-sided calibration criterion.  The formulation
with truncated divergences is convenient because it records explicitly the
passage to the principal value.

\begin{proposition}\label{prop:calibration}
Let $E\subset\R^n$ have locally finite fractional perimeter and let
$\zeta:\R^n\times\R^n\to[-1,1]$ be antisymmetric.  For $\rho>0$ set
\[
 \cD_{s,\rho}\zeta(x):=
 \int_{|x-y|>\rho}\zeta(x,y)|x-y|^{-n-s}\dd y.
\]
Assume that
\begin{enumerate}[label=\textup{(\roman*)}]
\item
\[
 \zeta(x,y)(\chi_E(y)-\chi_E(x))=|\chi_E(y)-\chi_E(x)|
\]
for a.e. pair $(x,y)$;
\item there exists $\cD_s\zeta\in L^1_{\rm loc}(E^c)$ such that
\[
 \cD_{s,\rho}\zeta\longrightarrow\cD_s\zeta
 \quad\hbox{in }L^1_{\rm loc}(E^c)\quad(\rho\downarrow0),
 \qquad
 \cD_s\zeta\le0\quad\hbox{a.e. in }E^c.
\]
\end{enumerate}
Then $E$ is outward minimizing for $P_s$.
\end{proposition}

\begin{proof}
Let $F\supset E$ with $F\setminus E\Subset\Omega$.  If
$P_s(F,\Omega)=+\infty$ there is nothing to prove, so assume it is finite and
put $u:=\chi_F-\chi_E$.  Then $u\ge0$ and $\operatorname{spt}u\Subset\Omega$.
For $\rho>0$ denote by $P_{s,\rho}$ the expression in
\eqref{eq:perimeter-Q} with the additional restriction $|x-y|>\rho$.
Since $|\zeta|\le1$, the saturation identity gives
\begin{align*}
 P_{s,\rho}(F,\Omega)-P_{s,\rho}(E,\Omega)
 &\ge \frac12\int_{\cQ_\Omega\cap\{|x-y|>\rho\}}
 \zeta(x,y)(u(y)-u(x))|x-y|^{-n-s}\dd x\dd y.
\end{align*}
Because $u$ vanishes outside a compact subset of $\Omega$, the last integral
is unchanged if $\cQ_\Omega$ is replaced by $\R^{2n}$.  It is absolutely
convergent for fixed $\rho$, and antisymmetry yields
\[
 \frac12\iint_{|x-y|>\rho}
 \zeta(x,y)(u(y)-u(x))|x-y|^{-n-s}\dd x\dd y
 =-\int_{E^c}u(x)\cD_{s,\rho}\zeta(x)\dd x.
\]
Letting $\rho\downarrow0$, monotone convergence gives
$P_{s,\rho}(A,\Omega)\uparrow P_s(A,\Omega)$ for $A=E,F$, while assumption
\textup{(ii)} gives convergence of the right-hand side.  Consequently
\[
 P_s(F,\Omega)-P_s(E,\Omega)
 \ge-\int_{E^c}u(x)\cD_s\zeta(x)\dd x\ge0.
\]
\end{proof}

\begin{corollary}\label{cor:foliation}
Let $\Psi:\R^n\to\R$ be smooth away from a closed null set and set
$E_a=\{\Psi\le a\}$.  Assume that $E_0$ has locally finite fractional
perimeter, that $\partial E_a$ is smooth for every $a>0$,
\[
 H_s^{E_a}>0\qquad\text{on }\partial E_a,
 \qquad a>0,
\]
and, for
\[
 \zeta(x,y):=\sgn(\Psi(x)-\Psi(y)),
\]
the truncated divergences $\cD_{s,\rho}\zeta$ converge in
$L^1_{\rm loc}(\{\Psi>0\})$.  Then $E_0=\{\Psi\le0\}$ is outward minimizing.
\end{corollary}

\begin{proof}
The saturation condition is immediate.  If $x\in E_0^c$ and
$a=\Psi(x)>0$, then, for a.e. $y$,
\[
 \zeta(x,y)=\chi_{E_a}(y)-\chi_{E_a^c}(y).
\]
Therefore
\[
 \cD_s\zeta(x)
 =\PV\int_{\R^n}\frac{\chi_{E_a}(y)-\chi_{E_a^c}(y)}{|x-y|^{n+s}}\dd y
 =-H_s^{E_a}(x)<0.
\]
The conclusion follows from \Cref{prop:calibration}.
\end{proof}

The final application of this criterion to the level sets of $\Phi$ and $-\Phi$
is postponed to \Cref{sec:final-summary}, after the strict barrier theorem has
been proved.

\section{Fractional mean curvature of the perturbed cone}\label{sec:curvature}

We now prove \Cref{thm:barrier}.  Since $E_1^+$ is invariant under
$O(4)\times O(4)$, every boundary point is equivalent to
\begin{equation}\label{eq:p-y}
 p(y):=(x(y)e_1,ye_1),
 \qquad x(y):=(1+y^4)^{1/4},
 \qquad y\ge0.
\end{equation}
The classical mean curvature with the sign convention \eqref{eq:Hs-def} is
\begin{equation}\label{eq:classical-H}
 \cH(y)=
 \frac{3\bigl(x(y)^2-y^2\bigr)^3}
 {\sqrt{x(y)^2+y^2}\,\bigl(x(y)^4-x(y)^2y^2+y^4\bigr)^{3/2}}.
\end{equation}
Hence
\begin{equation}\label{eq:H-positive-classical}
 \cH(y)>0\quad\forall y\ge0,
 \qquad
 \cH(y)=\frac{3}{8\sqrt2}y^{-13}+O(y^{-17})\quad(y\to\infty).
\end{equation}
Indeed, $x(y)^4-y^4=1$, and therefore
\[
 x(y)^2-y^2=\frac{x(y)^4-y^4}{x(y)^2+y^2}
 =\frac{1}{x(y)^2+y^2}>0.
\]
The denominator in \eqref{eq:classical-H} is strictly positive for every
$y\ge0$, including $y=0$.

\begin{lemma}\label{lem:compact-part}
For every $Y<\infty$ there is $s_Y\in(0,1)$ such that
\[
 H_s^{E_1^+}(p(y))>0
 \qquad\forall y\in[0,Y],\quad s\in(s_Y,1).
\]
\end{lemma}

\begin{proof}
On compact families of uniformly $C^2$ hypersurfaces,
\[
 (1-s)H_s^E(p)\to \kappa_n H_{\rm cl}^E(p)
\]
uniformly as $s\uparrow1$, with $\kappa_n>0$; see, for instance, \cite{AbatangeloValdinoci,CesaroniDeLucaNovagaPonsiglione}.  Since
\eqref{eq:H-positive-classical} gives a positive minimum of the classical mean
curvature on $[0,Y]$, the claim follows.
\end{proof}

It remains to handle $y\to\infty$.  Set
\begin{equation}\label{eq:eta-scaling}
 \eta:=y^{-4},
 \qquad
 E_\eta:=\{|X|^4\le |Y|^4+\eta\},
 \qquad
 P_\eta:=(a_\eta e_1,e_1),
 \qquad
 a_\eta:=(1+\eta)^{1/4}.
\end{equation}
Then
\begin{equation}\label{eq:H-scaling-eta}
 H_s^{E_1^+}(p(y))=y^{-s}H_s^{E_\eta}(P_\eta).
\end{equation}

\begin{proposition}\label{prop:uniform-expansion}
There exist $s_*\in(0,1)$, $\eta_0>0$, $c_*>0$ and $C_*<\infty$ such that, for
$s\in[s_*,1)$ and $0<\eta<\eta_0$,
\begin{equation}\label{eq:main-Heta-estimate}
 H_s^{E_\eta}(P_\eta)\ge A_s\eta-C_*\eta^2,
 \qquad
 A_s\ge c_*.
\end{equation}
In particular, after decreasing $\eta_0$ if necessary,
\[
 H_s^{E_\eta}(P_\eta)>0
 \qquad\forall s\in[s_*,1),\quad0<\eta<\eta_0.
\]
\end{proposition}

The proof occupies the rest of this section.

\subsection{Laplace--Gaussian representation}

Let $Z_t,W_t$ be independent Gaussian variables in $\R^4$ with density
\[
 \gamma_t(z)=\left(\frac{t}{\pi}\right)^2e^{-t|z|^2}.
\]
Define
\begin{equation}\label{eq:Pt-def}
 P_t(\eta):=\mathbb P\left(
 |a_\eta e_1+Z_t|^4-a_\eta^4
 > |e_1+W_t|^4-1
 \right).
\end{equation}

\begin{lemma}\label{lem:gaussian-representation}
For $0<s<1$ and $\eta\ge0$,
\begin{equation}\label{eq:Heta-Gaussian}
 H_s^{E_\eta}(P_\eta)=
 \frac{\pi^4}{\Gamma(4+s/2)}
 \int_0^\infty t^{-1+s/2}\bigl(2P_t(\eta)-1\bigr)\dd t.
\end{equation}
Moreover $P_t(0)=1/2$ for every $t>0$.
\end{lemma}

\begin{proof}
Let
\[
 \sigma_\eta(q):=
 \chi_{E_\eta^c}(P_\eta+q)-\chi_{E_\eta}(P_\eta+q).
\]
Choose the unit normal $\nu_\eta$ at $P_\eta$ with the orientation for which
$\sigma_\eta(q)=\sgn(q\cdot\nu_\eta)$ to first order, and put
$\sigma_\eta^0(q):=\sgn(q\cdot\nu_\eta)$.  Since
$\partial E_\eta$ is $C^2$ near $P_\eta$, there exist $r_0,C>0$ such that
\[
 \sigma_\eta(q)\ne\sigma_\eta^0(q),\quad |q|<r_0,
 \qquad\Longrightarrow\qquad
 |q\cdot\nu_\eta|\le C|q-(q\cdot\nu_\eta)\nu_\eta|^2.
\]
Consequently the set where the two signs differ occupies, in a spherical shell
of radius $r$, a proportion $O(r)$.  Hence
\[
 \int_{B_{r_0}}|\sigma_\eta(q)-\sigma_\eta^0(q)|
 |q|^{-8-s}\dd q<\infty
\]
precisely because $s<1$.  The same difference is integrable at infinity.
The principal value of the halfspace term is zero, and therefore
\begin{equation}\label{eq:H-absolute-difference}
 H_s^{E_\eta}(P_\eta)
 =\int_{\R^8}(\sigma_\eta(q)-\sigma_\eta^0(q))|q|^{-8-s}\dd q,
\end{equation}
where the right-hand side is absolutely convergent.

Using the Laplace identity
\[
 |q|^{-8-s}=\frac1{\Gamma(4+s/2)}
 \int_0^\infty t^{3+s/2}e^{-t|q|^2}\dd t
\]
and Fubini's theorem in \eqref{eq:H-absolute-difference}, justified by the absolute convergence of the iterated integral, we obtain
\[
 H_s^{E_\eta}(P_\eta)
 =\frac1{\Gamma(4+s/2)}\int_0^\infty t^{3+s/2}
 \int_{\R^8}(\sigma_\eta(q)-\sigma_\eta^0(q))e^{-t|q|^2}\dd q\dd t.
\]
The centered Gaussian integral of $\sigma_\eta^0$ vanishes by reflection in
$\nu_\eta^\perp$.  Since the density of $(Z_t,W_t)$ is
$(t/\pi)^4e^{-t|q|^2}$, the remaining inner integral equals
$\pi^4t^{-4}(2P_t(\eta)-1)$, proving \eqref{eq:Heta-Gaussian}.
When $\eta=0$, the random variables
$|e_1+Z_t|^4-1$ and $|e_1+W_t|^4-1$ are independent and have the same
continuous law, so $P_t(0)=1/2$.
\end{proof}

\begin{lemma}\label{lem:first-derivative}
For every $t>0$,
\begin{equation}\label{eq:Pprime-Bessel}
 \partial_\eta P_t(0)=
 \frac{e^{-t}}{16}\left(2tI_0(t)-(2t+1)I_1(t)\right),
\end{equation}
where $I_0$ and $I_1$ are modified Bessel functions of the first kind.
\end{lemma}

\begin{proof}
Write
\[
 U_a(z):=|ae_1+z|^4-a^4,
 \qquad a_\eta=(1+\eta)^{1/4}.
\]
For $\delta>0$, first restrict the coarea computation to the region
$|e_1+z|\ge\delta$ and then differentiate the resulting fixed-boundary
integral.  The omitted part is $O(\delta^4)$, uniformly for $\eta$ near zero,
while its differentiated integrand is $O(\delta^3)$ after the radial
Jacobian is included.  Letting $\delta\downarrow0$ therefore justifies
differentiation and gives
\begin{equation}\label{eq:Pprime-coarea}
 \partial_\eta P_t(0)
 =\int_{\R^4}\int_{\{U_1(z)=U_1(w)\}}
 \frac{G(z)}{|\nabla U_1(z)|}\gamma_t(z)\gamma_t(w)
 \dd\cH^3(z)\dd w,
\end{equation}
where
\[
 G(z):=\left.\partial_\eta U_{a_\eta}(z)\right|_{\eta=0}
 =|e_1+z|^2(1+z_1)-1.
\]
Indeed $a_\eta'(0)=1/4$ and
$\partial_aU_a(z)=4|ae_1+z|^2(a+z_1)-4a^3$.

Use polar coordinates centered at $-e_1$:
\[
 e_1+z=r\omega,
 \qquad e_1+w=r\nu,
 \qquad r>0,
 \qquad \omega,\nu\in\Sph^3.
\]
Then $U_1(z)=r^4-1$, $|\nabla U_1(z)|=4r^3$, and
\[
 \gamma_t(z)\gamma_t(w)=
 \left(\frac{t}{\pi}\right)^4e^{-2t(r^2+1)}
 e^{2tr(\omega_1+\nu_1)}.
\]
Moreover $1+z_1=r\omega_1$ and $|e_1+z|=r$, whence the corrected identity
\begin{equation}\label{eq:G-correct}
 G(z)=r^3\omega_1-1.
\end{equation}
Thus
\begin{equation}\label{eq:Pprime-spherical}
 \partial_\eta P_t(0)=
 \frac14\left(\frac{t}{\pi}\right)^4
 \int_0^\infty r^3e^{-2t(r^2+1)}
 \bigl(r^3A_1(2tr)-A_0(2tr)\bigr)A_0(2tr)\dd r,
\end{equation}
where
\[
 A_0(\rho):=\int_{\Sph^3}e^{\rho\omega_1}\dd\sigma(\omega),
 \qquad A_1(\rho):=A_0'(\rho).
\]
The standard integral representation of the modified Bessel functions, see \cite[Section~10.32]{NIST}, gives
\[
 A_0(\rho)=\frac{4\pi^2I_1(\rho)}{\rho},
 \qquad
 A_1(\rho)=4\pi^2\left(\frac{I_0(\rho)}{\rho}
 -\frac{2I_1(\rho)}{\rho^2}\right).
\]
Substitution in \eqref{eq:Pprime-spherical} gives, with all Bessel functions
in the following line evaluated at $2tr$,
\begin{equation}\label{eq:Pprime-radial-bessel}
 \partial_\eta P_t(0)=e^{-2t}\int_0^\infty e^{-2tr^2}
 \left(t^2r^4I_0I_1-tr^3I_1^2-t^2rI_1^2\right)\dd r.
\end{equation}
To evaluate the integral, set
\[
 W(\alpha,\beta):=
 \int_0^\infty r e^{-\alpha r^2}I_1(\beta r)^2\dd r
 =\frac1{2\alpha}e^{\beta^2/(2\alpha)}
 I_1\!\left(\frac{\beta^2}{2\alpha}\right).
\]
Using $I_1'(x)=I_0(x)-I_1(x)/x$, differentiation under the integral gives
\[
 \int_0^\infty e^{-\alpha r^2}
 \left(t^2r^4I_0(\beta r)I_1(\beta r)
 -tr^3I_1(\beta r)^2-t^2rI_1(\beta r)^2\right)\dd r
\]
\[
 =-\frac{t^2}{2}W_{\alpha\beta}(\alpha,\beta)
 +\frac t2W_\alpha(\alpha,\beta)-t^2W(\alpha,\beta)
 \quad\text{when }\beta=2t.
\]
Evaluating at $\alpha=\beta=2t$ and differentiating the explicit Weber
formula yields
\[
 -\frac{t^2}{2}W_{\alpha\beta}+\frac t2W_\alpha-t^2W
 =\frac{e^t}{16}\bigl(2tI_0(t)-(2t+1)I_1(t)\bigr).
\]
Together with \eqref{eq:Pprime-radial-bessel}, this proves
\eqref{eq:Pprime-Bessel}.  All differentiations are justified by the Gaussian
factor and standard local bounds for $I_0,I_1$.
\end{proof}

\begin{lemma}\label{lem:Bessel-positive}
For every $t>0$,
\begin{equation}\label{eq:Bessel-positive}
 2tI_0(t)-(2t+1)I_1(t)>0.
\end{equation}
Consequently
\begin{equation}\label{eq:As}
 A_s:=\frac{\pi^4}{8\Gamma(4+s/2)}
 \int_0^\infty t^{-1+s/2}e^{-t}
 \left(2tI_0(t)-(2t+1)I_1(t)\right)\dd t
\end{equation}
is continuous and strictly positive on every compact subinterval of $(0,1]$.
In particular, for $s_*<1$ sufficiently close to $1$, there exists $c_*>0$ such
that
\[
 A_s\ge c_*\qquad\forall s\in[s_*,1].
\]
\end{lemma}

\begin{proof}
Let $R(t):=I_1(t)/I_0(t)$ and $b(t):=2t/(2t+1)$.  The identities
$I_0'=I_1$ and $I_1'=I_0-I_1/t$ give
\[
 R'=1-R^2-\frac Rt.
\]
Since $R(t)=t/2+O(t^3)$ while $b(t)=2t+O(t^2)$, $R<b$ for small positive $t$.
If $R-b$ had a first zero $t_0>0$, then $(R-b)'(t_0)\ge0$.  But at $R=b$,
\[
 (R-b)'=1-b^2-\frac bt-b'=-\frac{3}{(2t+1)^2}<0,
\]
a contradiction.  Hence $R<b$ for all $t>0$, which is exactly
\eqref{eq:Bessel-positive}.

At zero,
\[
 e^{-t}\bigl(2tI_0(t)-(2t+1)I_1(t)\bigr)=\frac32t+O(t^2).
\]
At infinity, the standard large-argument expansions of $I_0$ and $I_1$ \cite[Section~10.40]{NIST} give
\[
 e^{-t}\bigl(2tI_0(t)-(2t+1)I_1(t)\bigr)=O(t^{-3/2}).
\]
These two bounds imply uniform integrability of the defining integral for $A_s$
on compact intervals of $s\le1$, and continuity follows by dominated convergence.
The integrand is positive, hence $A_s>0$.
\end{proof}

\subsection{A one-sided Taylor estimate}\label{subsec:one-sided}

The bounded and the large Gaussian scales require different coordinates.  At
bounded scales we use the exact one-dimensional density of a noncentral
chi-square variable.  At large scales we use a normal graph, but keep the
Taylor remainder in $\eta$ throughout the argument.

\begin{lemma}\label{lem:bounded-scales}
There exist $\eta_0>0$ and $C<\infty$ such that
\begin{equation}\label{eq:small-t-second}
 |\partial_\eta^2P_t(\eta)|\le Ct,
 \qquad 0<t\le1,
 \quad 0\le\eta\le\eta_0.
\end{equation}
In particular,
\[
 2P_t(\eta)-1
 \ge2\eta\partial_\eta P_t(0)-C\eta^2t.
\]
\end{lemma}

\begin{proof}
Let $\Xi,\Omega\in\R^4$ be independent with density
$\pi^{-2}e^{-|\xi|^2}$.  After writing $Z_t=t^{-1/2}\Xi$ and
$W_t=t^{-1/2}\Omega$, set
\[
 X_\eta:=|\sqrt t\,a_\eta e_1+\Xi|^2,
 \qquad
 Y:=|\sqrt t\,e_1+\Omega|^2.
\]
The event in \eqref{eq:Pt-def} is exactly
\[
 X_\eta^2>Y^2+t^2\eta,
 \qquad\text{or equivalently}\qquad
 X_\eta>\sqrt{Y^2+t^2\eta}.
\]
If $\lambda\ge0$, the variable $|\sqrt\lambda e_1+\Xi|^2$ has density
\begin{equation}\label{eq:noncentral-density}
 f_\lambda(x)
 =e^{-(x+\lambda)}\sqrt{\frac{x}{\lambda}}
 I_1(2\sqrt{\lambda x})
 =x h_\lambda(x),
 \qquad x>0,
\end{equation}
with the value at $\lambda=0$ understood by continuity, and
\begin{equation}\label{eq:h-series}
 h_\lambda(x)=e^{-(x+\lambda)}
 \sum_{k=0}^\infty\frac{(\lambda x)^k}{k!(k+1)!}.
\end{equation}
Let
\[
 Q(\lambda,x):=\int_x^\infty f_\lambda(r)\dd r,
 \qquad
 \lambda_\eta:=ta_\eta^2=t\sqrt{1+\eta},
 \qquad
 x_\eta(y):=\sqrt{y^2+t^2\eta}.
\]
Conditioning on $Y$ gives the fixed-domain formula
\begin{equation}\label{eq:P-radial-fixed-domain}
 P_t(\eta)=\int_0^\infty
 Q(\lambda_\eta,x_\eta(y))f_t(y)\dd y.
\end{equation}
Formula \eqref{eq:h-series}, or equivalently the standard bounds for $I_1$,
shows that for $0\le\lambda\le2$ the functions $h_\lambda$ and their
$\lambda$-derivatives up to order two and $x$-derivatives up to order one are
bounded by $C(1+x)^Me^{-cx}$.  The same type of bound holds for
$Q_\lambda$ and $Q_{\lambda\lambda}$.  For $\eta>0$ these bounds justify two differentiations in
\eqref{eq:P-radial-fixed-domain}.  To include $\eta=0$, differentiate first
on $[\delta,\infty)$ in the $y$-variable, combine the two terms containing
$x_\eta''$ before taking the limit, and then let $\delta\downarrow0$.  The
combined expression below is dominated by an integrable function independent
of $\delta$ and $\eta$.  Thus the second derivative extends continuously from
$\eta>0$ to the right endpoint $\eta=0$.

We record the cancellation at $y=0$, which is the point missed by an argument
based on a $C^2$ density for the difference variable.  Since
$Q_x=-xh_\lambda(x)$ and $Q_{xx}=-h_\lambda(x)-xh_\lambda'(x)$,
\[
 Q_{xx}(x_\eta')^2+Q_x x_\eta''
 =-\frac{t^4}{4x_\eta}h_{\lambda_\eta}'(x_\eta),
\]
where
\[
 x_\eta'=\frac{t^2}{2x_\eta},
 \qquad
 x_\eta''=-\frac{t^4}{4x_\eta^3}.
\]
Similarly,
\[
 2Q_{\lambda x}\lambda_\eta'x_\eta'
 =-\lambda_\eta't^2
 \partial_\lambda h_{\lambda_\eta}(x_\eta).
\]
Thus
\begin{align}\label{eq:Psecond-exact-bounded}
 \partial_\eta^2P_t(\eta)
 =\int_0^\infty &\bigg[
 Q_{\lambda\lambda}(\lambda_\eta,x_\eta)(\lambda_\eta')^2
 +Q_\lambda(\lambda_\eta,x_\eta)\lambda_\eta'' \\
 &\quad-\lambda_\eta't^2
 \partial_\lambda h_{\lambda_\eta}(x_\eta)
 -\frac{t^4}{4x_\eta}h_{\lambda_\eta}'(x_\eta)
 \bigg]f_t(y)\dd y.
\end{align}
Here and below $x_\eta=x_\eta(y)$.  Since $x_\eta\ge y$ and
$f_t(y)=yh_t(y)$, the last term is integrable uniformly even when
$\eta=0$.  Moreover
$|\lambda_\eta'|\le Ct$ and $|\lambda_\eta''|\le Ct$.  Applying the bounds
following \eqref{eq:h-series} in \eqref{eq:Psecond-exact-bounded} gives
\[
 |\partial_\eta^2P_t(\eta)|
 \le C(t+t^2+t^3+t^4)\le Ct.
\]
Taylor's formula and $P_t(0)=1/2$ finish the proof.
\end{proof}

\begin{lemma}\label{lem:large-scales}
There exist $\eta_0>0$, $\eps_0>0$ and $C<\infty$ such that, for
$0<\eps\le\eps_0$ and $0\le\eta\le\eta_0$, one can write
\begin{equation}\label{eq:large-scale-decomposition}
 P_{\eps^{-2}}(\eta)-\frac12
 =B_\eps(\eta)+\eps^3\Psi_\eps(\eta),
\end{equation}
where
\begin{equation}\label{eq:B-eps-classical}
 B_\eps(\eta)=\frac{\eps}{4\sqrt\pi}
 \frac{3(a_\eta^2-1)^3}
 {\sqrt{a_\eta^2+1}\,(a_\eta^4-a_\eta^2+1)^{3/2}}
 \ge0,
\end{equation}
$B_\eps(0)=B_\eps'(0)=0$, and
\begin{equation}\label{eq:Psi-second-bound}
 \sup_{0<\eps\le\eps_0}\sup_{0\le\eta\le\eta_0}
 |\Psi_\eps''(\eta)|\le C.
\end{equation}
Consequently
\begin{equation}\label{eq:large-t-Taylor}
 2P_t(\eta)-1
 \ge2\eta\partial_\eta P_t(0)-C\eta^2t^{-3/2}
 \qquad(t\ge\eps_0^{-2}).
\end{equation}
\end{lemma}

\begin{proof}
Put $z=\eps\xi$, $w=\eps\omega$ and
$F_{\eta,\eps}(\xi,\omega):=F_\eta(\eps\xi,\eps\omega)$.  Then
\[
 P_{\eps^{-2}}(\eta)=\pi^{-4}\int_{\R^8}
 \1_{\{F_{\eta,\eps}>0\}}e^{-|\xi|^2-|\omega|^2}\dd\xi\dd\omega,
\]
and the following expansion is exact:
\begin{align}\label{eq:F-expansion}
F_{\eta,\eps}
&=4\eps(a_\eta^3\xi_1-\omega_1)
+\eps^2\{2a_\eta^2|\xi|^2+4a_\eta^2\xi_1^2
       -2|\omega|^2-4\omega_1^2\} \\
&\quad+4\eps^3\{a_\eta\xi_1|\xi|^2-\omega_1|\omega|^2\}
+\eps^4\{|\xi|^4-|\omega|^4\}.\nonumber
\end{align}
Introduce the orthogonal coordinates
\[
 N_\eta=(a_\eta^6+1)^{1/2},\qquad
 \alpha=\frac{a_\eta^3\xi_1-\omega_1}{N_\eta},\qquad
 \beta=\frac{\xi_1+a_\eta^3\omega_1}{N_\eta},
\]
so that
\[
 \xi_1=\frac{a_\eta^3\alpha+\beta}{N_\eta},
 \qquad
 \omega_1=\frac{-\alpha+a_\eta^3\beta}{N_\eta}.
\]
Write $\tau=(\beta,\xi',\omega')\in\R^7$, $u=|\xi'|^2$ and
$v=|\omega'|^2$.  In these variables
\begin{equation}\label{eq:adapted-F-expansion}
 F_{\eta,\eps}=4\eps N_\eta\alpha
 +\eps^2R_{2,\eta}(\alpha,\tau)
 +\eps^3R_{3,\eta}(\alpha,\tau)
 +\eps^4R_{4,\eta}(\alpha,\tau),
\end{equation}
where the $R_{j,\eta}$ are polynomials with coefficients smooth in $\eta$.
At $\alpha=0$,
\begin{equation}\label{eq:R2-adapted}
 R_{2,\eta}(0,\tau)
 =2a_\eta^2u-2v-\frac{6a_\eta^2\eta}{N_\eta^2}\beta^2.
\end{equation}

We give the localization details needed for uniformity.  Set
$L_\eps:=\eps^{-1/8}$.  On
$\{|\alpha|+|\tau|\le2L_\eps\}$, differentiation of
\eqref{eq:adapted-F-expansion} gives
\[
 \partial_\alpha F_{\eta,\eps}
 =4\eps N_\eta+O(\eps^2(1+L_\eps)^3)
 \ge2\eps N_\eta
\]
for $\eps_0$ sufficiently small.  The signs at $\alpha=\pm L_\eps$ are
opposite, since the linear term has size $\eps L_\eps$ whereas all remaining
terms are $o(\eps L_\eps)$.  Hence the zero set in this cylinder is a unique
normal graph
\[
 \alpha=h_{\eta,\eps}(\tau).
\]
Write $h_{\eta,\eps}=\eps k_{\eta,\eps}$ in the implicit equation
$F_{\eta,\eps}(h_{\eta,\eps},\tau)=0$.  Since the derivative with respect to
$k$ is $4N_\eta+O(\eps(1+|\tau|)^3)$, the quantitative implicit-function
theorem applies uniformly on the cylinder.  Expanding the equation through
order $\eps^2$ and differentiating it twice with respect to $\eta$ gives
\begin{equation}\label{eq:h-expansion-rigorous}
 h_{\eta,\eps}(\tau)
 =\eps H_{1,\eta}(\tau)+\eps^2H_{2,\eta}(\tau)
 +\eps^3K_{\eta,\eps}(\tau),
\end{equation}
where, for $j=0,1,2$,
\begin{equation}\label{eq:h-eta-bounds}
 |\partial_\eta^jH_{1,\eta}|\le C(1+|\tau|)^2,
 \quad
 |\partial_\eta^jH_{2,\eta}|\le C(1+|\tau|)^3,
 \quad
 |\partial_\eta^jK_{\eta,\eps}|\le C(1+|\tau|)^M.
\end{equation}
The first coefficient is
\begin{equation}\label{eq:H1-explicit}
 H_{1,\eta}
 =\frac{v-a_\eta^2u}{2N_\eta}
 +\frac{3a_\eta^2\eta}{2N_\eta^3}\beta^2.
\end{equation}
The coefficient equation at order $\eps^2$ is
\begin{equation}\label{eq:H2-explicit-recursion}
 H_{2,\eta}=-\frac{\partial_\alpha R_{2,\eta}(0,\tau)H_{1,\eta}
                    +R_{3,\eta}(0,\tau)}{4N_\eta}.
\end{equation}
Furthermore, $R_{2,\eta}(0,\tau)$ and $H_{1,\eta}$ are even in $\beta$,
whereas $\partial_\alpha R_{2,\eta}(0,\tau)$ and
$R_{3,\eta}(0,\tau)$ are odd in $\beta$.  Formula
\eqref{eq:H2-explicit-recursion} therefore gives
\begin{equation}\label{eq:H2-odd}
 H_{2,\eta}(-\beta,\xi',\omega')
 =-H_{2,\eta}(\beta,\xi',\omega').
\end{equation}

For completeness, we also record how the Gaussian tails are handled.  A
smooth cutoff equal to one on $|\alpha|+|\tau|\le L_\eps$ and zero outside
$|\alpha|+|\tau|\le2L_\eps$ is inserted.  On its support the preceding graph
description is valid.  The gradient in the original scaled variables is
\[
 \nabla_\xi F_{\eta,\eps}
 =4\eps|a_\eta e_1+\eps\xi|^2(a_\eta e_1+\eps\xi),
 \qquad
 \nabla_\omega F_{\eta,\eps}
 =-4\eps|e_1+\eps\omega|^2(e_1+\eps\omega).
\]
Thus the only common critical point is
$(-a_\eta e_1/\eps,-e_1/\eps)$, and its value is $-\eta$.  It does not
belong to the zero set when $\eta>0$; when $\eta=0$ it lies at Gaussian
distance $O(\eps^{-1})$.  Away from a fixed neighborhood of that point, the
coarea formulas for the first two $\eta$-derivatives have a nonvanishing
gradient denominator and produce only polynomial factors.  In the remaining
neighborhood use
\[
 p=a_\eta e_1+\eps\xi,
 \qquad q=e_1+\eps\omega,
 \qquad u=|p|^4,
 \qquad v=|q|^4.
\]
Since in four dimensions $r^3\dd r=\frac14\dd(r^4)$, the radial Jacobian in
$(u,v)$ is constant and the event is the half-plane $u-v>\eta$.  After the
angular integrations, each radial density has an expansion
$g(u)=g(0)+O(u^{1/2})$ and satisfies $g'(u)=O(u^{-1/2})$ at the origin, with
all coefficients bounded by a polynomial in $\eps^{-1}$ times
$e^{-c/\eps^2}$.  For $d\ge0$ the density of $u-v$ is
$k(d)=\int_0^\infty g_1(v+d)g_2(v)\dd v$; hence
$k'(0+)=\int_0^\infty g_1'(v)g_2(v)\dd v$ is finite.  This gives the two
right $\eta$-derivatives required at $\eta=0$, with the same exponentially
small bound.  Combining the two regions, the complementary contribution
$\mathcal R_\eps(\eta)$ is $C^2$ and satisfies
\begin{equation}\label{eq:gaussian-tail-C2}
 \sup_{0\le\eta\le\eta_0}|\mathcal R_\eps''(\eta)|
 \le Ce^{-cL_\eps^2}.
\end{equation}
Possible additional roots of the quartic in $\alpha$ lie outside the cutoff
cylinder; they are included in $\mathcal R_\eps$.  This is the reason for
localizing before replacing the superlevel set by a single half-line.
Since $e^{-cL_\eps^2}=e^{-c\eps^{-1/4}}=O(\eps^N)$ for every $N$, the
remainder in \eqref{eq:gaussian-tail-C2} may be absorbed at order $\eps^3$
without losing two $\eta$-derivatives.

Let
\[
 q(r):=\pi^{-1/2}\int_r^\infty e^{-\alpha^2}\dd\alpha,
 \qquad
 \dd\gamma_7(\tau):=\pi^{-7/2}e^{-|\tau|^2}\dd\tau.
\]
After extending the localized graph smoothly outside the central region, the
preceding discussion gives
\begin{equation}\label{eq:P-via-threshold-rigorous}
 P_{\eps^{-2}}(\eta)
 =\int_{\R^7}q(h_{\eta,\eps}(\tau))\dd\gamma_7(\tau)
 +\mathcal R_\eps(\eta).
\end{equation}
On the central region one has $|h_{\eta,\eps}|\le C\eps^{3/4}$ for $\eps_0$ sufficiently small.  Hence
\[
 q(r)=\frac12-\frac r{\sqrt\pi}+r^3\rho(r),
\]
where $\rho$ is smooth.  Substitution of
\eqref{eq:h-expansion-rigorous} into
\eqref{eq:P-via-threshold-rigorous}, together with
\eqref{eq:h-eta-bounds} and the finiteness of all Gaussian moments, yields
\begin{equation}\label{eq:P-normal-form-decomposition}
 P_{\eps^{-2}}(\eta)-\frac12
 =-\frac\eps{\sqrt\pi}\int H_{1,\eta}\dd\gamma_7
 -\frac{\eps^2}{\sqrt\pi}\int H_{2,\eta}\dd\gamma_7
 +\eps^3\Psi_\eps(\eta),
\end{equation}
with \eqref{eq:Psi-second-bound}.  Notice that the $C^2$ estimate for
$\Psi_\eps$ follows from the $\eta$-derivative bounds in
\eqref{eq:h-eta-bounds}; a bare $O(\eps^5)$ remainder would not suffice for
this conclusion.

The second integral in \eqref{eq:P-normal-form-decomposition} vanishes by
\eqref{eq:H2-odd}.  From \eqref{eq:H1-explicit} and the Gaussian moments
$\int\beta^2\dd\gamma_7=1/2$ and
$\int u\dd\gamma_7=\int v\dd\gamma_7=3/2$, one obtains
\begin{equation}\label{eq:H1-average-classical}
 \int H_{1,\eta}\dd\gamma_7
 =-\frac14\,
 \frac{3(a_\eta^2-1)^3}
 {\sqrt{a_\eta^2+1}\,(a_\eta^4-a_\eta^2+1)^{3/2}}.
\end{equation}
Equations \eqref{eq:P-normal-form-decomposition} and
\eqref{eq:H1-average-classical} prove
\eqref{eq:large-scale-decomposition}--\eqref{eq:B-eps-classical}.
Finally Taylor's theorem, $B_\eps\ge0$, and
$B_\eps(0)=B_\eps'(0)=0$ give
\[
 P_{\eps^{-2}}(\eta)-P_{\eps^{-2}}(0)
 -\eta\partial_\eta P_{\eps^{-2}}(0)
 \ge-C\eta^2\eps^3,
\]
which is \eqref{eq:large-t-Taylor}.
\end{proof}

\begin{lemma}\label{lem:one-sided-gaussian}
There are $\eta_0>0$ and $C<\infty$ such that, for every
$0<\eta<\eta_0$ and every $t>0$,
\begin{equation}\label{eq:one-sided-Gaussian}
 2P_t(\eta)-1
 \ge2\eta\partial_\eta P_t(0)-C\eta^2m(t),
\end{equation}
where
\begin{equation}\label{eq:m-def}
 m(t)=
 \begin{cases}
 t, &0<t\le1,\\
 t^{-3/2}, &t\ge1.
 \end{cases}
\end{equation}
\end{lemma}

\begin{proof}
The range $0<t\le1$ is \Cref{lem:bounded-scales}.  The large-scale lemma
covers $t\ge\eps_0^{-2}$.  On the remaining compact interval
$1\le t\le\eps_0^{-2}$, the fixed-domain representation
\eqref{eq:P-radial-fixed-domain} and the same differentiation used in
\eqref{eq:Psecond-exact-bounded} give
$|\partial_\eta^2P_t(\eta)|\le C$; here the density parameters remain in a fixed compact interval depending only on $\eps_0$.  Since $t^{-3/2}$ is bounded below by a
positive constant on this compact interval, Taylor's theorem gives
\eqref{eq:one-sided-Gaussian} after increasing $C$.
\end{proof}


\subsection{Proof of the uniform expansion}

\begin{proof}[Proof of \Cref{prop:uniform-expansion}]
By \Cref{lem:gaussian-representation,lem:one-sided-gaussian},
\[
\begin{aligned}
H_s^{E_\eta}(P_\eta)
&\ge \frac{\pi^4}{\Gamma(4+s/2)}
 \int_0^\infty t^{-1+s/2}
 \left(2\eta\partial_\eta P_t(0)-C\eta^2m(t)\right)\dd t \\
&= A_s\eta
 -C\eta^2\frac{\pi^4}{\Gamma(4+s/2)}
 \int_0^\infty t^{-1+s/2}m(t)\dd t.
\end{aligned}
\]
The last integral is uniformly bounded for $s\in[s_*,1]$, because
\[
 \int_0^1 t^{-1+s/2}t\dd t=\frac{1}{1+s/2},
\]
and
\[
 \int_1^\infty t^{-1+s/2}t^{-3/2}\dd t
 =\frac{1}{3/2-s/2}.
\]
The factor $\Gamma(4+s/2)^{-1}$ is uniformly bounded on $[s_*,1]$.  Hence
\[
 H_s^{E_\eta}(P_\eta)\ge A_s\eta-C_*\eta^2.
\]
By \Cref{lem:Bessel-positive}, choosing $s_*<1$ sufficiently close to $1$ gives
$A_s\ge c_*>0$.  Finally decrease $\eta_0$ so that
$C_*\eta_0\le c_*/2$.
\end{proof}

\begin{corollary}\label{cor:infinity}
There exist $Y<\infty$, $s_*<1$ and $c>0$ such that
\[
 H_s^{E_1^+}(p(y))\ge c y^{-s-4}>0
 \qquad\forall y\ge Y,
 \quad s\in[s_*,1).
\]
\end{corollary}

\begin{proof}
Take $\eta=y^{-4}$ in \eqref{eq:H-scaling-eta}.  For $y$ large enough,
\Cref{prop:uniform-expansion} gives
\[
 H_s^{E_1^+}(p(y))=y^{-s}H_s^{E_\eta}(P_\eta)
 \ge y^{-s}\frac{c_*}{2}\eta
 =\frac{c_*}{2}y^{-s-4}.
\]
\end{proof}

\begin{proof}[Proof of \Cref{thm:barrier}]
Use \Cref{cor:infinity} for $y\ge Y$ and \Cref{lem:compact-part} on $[0,Y]$.
The statement for $E_1^-$ follows from the isometry $(X,Y)\mapsto(Y,X)$.
\end{proof}

\section{Proof of the main result}\label{sec:final-summary}

We now turn the strict barriers into two one-sided calibrations and conclude by submodularity.

\begin{proposition}\label{prop:summary-barriers}
There exists $s_b\in(0,1)$ such that, for every $s\in(s_b,1)$ and every
$a>0$,
\begin{equation}\label{eq:summary-positive-plus}
 H_s^{E_a^+}(p)>0
 \qquad\text{for every }p\in\partial E_a^+,
 \qquad
 E_a^+:=\{\Phi\le a\},
\end{equation}
and
\begin{equation}\label{eq:summary-positive-minus}
 H_s^{E_a^-}(p)>0
 \qquad\text{for every }p\in\partial E_a^-,
 \qquad
 E_a^-:=\{-\Phi\le a\}.
\end{equation}
\end{proposition}

\begin{proof}
By Corollary~\ref{cor:infinity} there exist $Y<\infty$, $s_\infty\in(0,1)$ and
$c_\infty>0$ such that
\[
 H_s^{E_1^+}(p(y))\ge c_\infty y^{-s-4}>0
 \qquad\forall y\ge Y,\quad s\in[s_\infty,1).
\]
For the same $Y$, Lemma~\ref{lem:compact-part} gives $s_Y\in(0,1)$ such that
\[
 H_s^{E_1^+}(p(y))>0
 \qquad\forall y\in[0,Y],\quad s\in(s_Y,1).
\]
Set
\[
 s_b:=\max\{s_\infty,s_Y\}.
\]
Since every point of $\partial E_1^+$ is, by the $O(4)\times O(4)$ symmetry, of
the form $p(y)$ for some $y\ge0$, the two estimates imply
\[
 H_s^{E_1^+}(p)>0
 \qquad\forall p\in\partial E_1^+,\quad s\in(s_b,1).
\]

Now fix $a>0$.  Because $\Phi$ is homogeneous of degree $4$,
\[
 E_a^+=a^{1/4}E_1^+.
\]
If $p\in\partial E_a^+$ and $q=a^{-1/4}p$, then $q\in\partial E_1^+$.  The scaling
identity in Lemma~\ref{lem:scaling} gives
\[
 H_s^{E_a^+}(p)
 =H_s^{a^{1/4}E_1^+}(a^{1/4}q)
 =a^{-s/4}H_s^{E_1^+}(q)>0.
\]
This proves \eqref{eq:summary-positive-plus}.  Finally, the exchange map
\[
 T(X,Y):=(Y,X)
\]
is an isometry of $\R^8$ and satisfies $T(E_a^+)=E_a^-$.  Fractional mean
curvature is invariant under isometries, hence \eqref{eq:summary-positive-minus}
follows from \eqref{eq:summary-positive-plus}.
\end{proof}

\begin{lemma}\label{lem:homogeneous-calibration}
Let
\[
 \zeta(x,y):=\sgn\bigl(\Phi(x)-\Phi(y)\bigr).
\]
Then, for every $0<s<1$, the truncated divergences
$\cD_{s,\rho}\zeta$ converge in
$L^1_{\rm loc}(C^c)$ as $\rho\downarrow0$.  Their limit satisfies
\begin{equation}\label{eq:Dzeta-leaf-curvature}
 \cD_s\zeta(x)=-H_s^{E_{\Phi(x)}^+}(x),
 \qquad x\in C^c,
\end{equation}
and
\begin{equation}\label{eq:Dzeta-homogeneous-bound}
 |\cD_{s,\rho}\zeta(x)|+|\cD_s\zeta(x)|
 \le C_s|x|^{-s}
\end{equation}
locally uniformly with respect to $\rho>0$.
\end{lemma}

\begin{proof}
For $x\in C^c$ put $a=\Phi(x)>0$.  For a.e. $y$,
\[
 \zeta(x,y)=\chi_{E_a^+}(y)-\chi_{(E_a^+)^c}(y),
\]
which gives \eqref{eq:Dzeta-leaf-curvature} once the principal value is shown
to exist.  We verify this uniformly after normalization.

The kernel is homogeneous and $\Phi$ has degree four, hence
\[
 \zeta(\lambda x,\lambda y)=\zeta(x,y),
 \qquad
 \cD_{s,\rho}\zeta(\lambda x)
 =\lambda^{-s}\cD_{s,\rho/\lambda}\zeta(x).
\]
Consider $x=\theta\in\Sph^7\cap C^c$.  The level hypersurface
$\{\Phi=\Phi(\theta)\}$ through $\theta$ has uniformly bounded $C^2$ geometry
as $\theta$ ranges in the closure of $\Sph^7\cap C^c$: on the limiting cone
$\{\Phi=0\}\cap\Sph^7$ one still has $\nabla\Phi\ne0$.  Subtracting the sign
of the tangent halfspace as in the proof of
\Cref{lem:gaussian-representation}, the contribution from a ball of radius
$r$ is bounded uniformly by
\[
 C\int_0^r\rho^{-s}\dd\rho=C r^{1-s}.
\]
The contribution outside a fixed ball is uniformly bounded by the integrable
tail of $|x-y|^{-8-s}$.  It follows that
$\cD_{s,\delta}\zeta(\theta)$ converges uniformly in
$\theta\in\Sph^7\cap C^c$ and is uniformly bounded there.

Writing $x=r\theta$ gives \eqref{eq:Dzeta-homogeneous-bound}.  Since
$|x|^{-s}\in L^1_{\rm loc}(\R^8)$, pointwise convergence away from the origin
and dominated convergence imply convergence in $L^1_{\rm loc}(C^c)$.
\end{proof}

\begin{proposition}\label{prop:summary-exterior}
For every $s\in(s_b,1)$, the cone
\[
 C=\{\Phi\le0\}
\]
is outward minimizing for $P_s$.
\end{proposition}

\begin{proof}
Apply Corollary~\ref{cor:foliation} with
\[
 \Psi=\Phi,
 \qquad
 E_a=\{\Psi\le a\}=E_a^+.
\]
For $a>0$, the boundary $\partial E_a^+$ is smooth because
$\nabla\Phi\ne0$ on $\{\Phi=a\}$.  Indeed, $\nabla\Phi(X,Y)=0$ only when
$X=Y=0$, whereas $0\notin\{\Phi=a\}$ for $a>0$.  The cone has locally finite classical perimeter; the standard translation estimate for characteristic functions therefore implies locally finite fractional perimeter for every $0<s<1$.  The positivity assumption in
Corollary~\ref{cor:foliation} is precisely \eqref{eq:summary-positive-plus},
and the required $L^1_{\rm loc}$ convergence of the truncated divergences is
\Cref{lem:homogeneous-calibration}.  Therefore
$E_0=\{\Phi\le0\}=C$ is outward minimizing.

\end{proof}

\begin{proposition}\label{prop:summary-interior}
For every $s\in(s_b,1)$, the cone $C$ is inward minimizing for $P_s$.
\end{proposition}

\begin{proof}
We apply the previous argument to $-\Phi$.
The positive leaves are
\[
 E_a^-=\{-\Phi\le a\}=\{|Y|^4\le |X|^4+a\},
 \qquad a>0,
\]
and their strict mean-convexity is \eqref{eq:summary-positive-minus}.  Hence
\[
 \{-\Phi\le0\}=\{\Phi\ge0\}
\]
is outward minimizing.  Since $\{\Phi\ge0\}$ differs from $C^c=\{\Phi>0\}$ only
by the boundary $\{\Phi=0\}$, a Lebesgue null set, the same outward minimality
holds for $C^c$.

Let now $F\subset C$ and $C\setminus F\Subset\Omega$.  Then
\[
 F^c\supset C^c,
 \qquad
 F^c\setminus C^c=C\setminus F\Subset\Omega.
\]
Outward minimality of $C^c$ gives
\[
 P_s(C^c,\Omega)\le P_s(F^c,\Omega).
\]
Using complement invariance from Lemma~\ref{lem:lattice},
\[
 P_s(C,\Omega)=P_s(C^c,\Omega)
 \le P_s(F^c,\Omega)=P_s(F,\Omega).
\]
This is exactly inward minimality of $C$.
\end{proof}

\begin{theorem}[Minimality of the Simons cone]\label{thm:summary-main}
Let $s_0:=s_b$.  For every $s\in(s_0,1)$, every bounded open set
$\Omega\subset\R^8$, and every measurable competitor $F\subset\R^8$ satisfying
\[
 F\Delta C\Subset\Omega,
\]
one has
\[
 P_s(C,\Omega)\le P_s(F,\Omega).
\]
In particular, Theorem~\ref{thm:main} follows.
\end{theorem}

\begin{proof}
Define the two one-sided competitors
\[
 F^+:=F\cup C,
 \qquad
 F^-:=F\cap C.
\]
Then
\[
 F^+\supset C,
 \qquad
 F^+\setminus C=F\setminus C\subset F\Delta C\Subset\Omega,
\]
and
\[
 F^-\subset C,
 \qquad
 C\setminus F^-=C\setminus F\subset F\Delta C\Subset\Omega.
\]
By Propositions~\ref{prop:summary-exterior} and~\ref{prop:summary-interior},
\[
 P_s(C,\Omega)\le P_s(F^+,\Omega),
 \qquad
 P_s(C,\Omega)\le P_s(F^-,\Omega).
\]
Adding the two inequalities yields
\begin{equation}\label{eq:summary-two-sided}
 2P_s(C,\Omega)
 \le
 P_s(F\cup C,\Omega)+P_s(F\cap C,\Omega).
\end{equation}
The submodularity inequality \eqref{eq:submod}, applied with $A=F$ and $B=C$,
gives
\begin{equation}\label{eq:summary-submod}
 P_s(F\cup C,\Omega)+P_s(F\cap C,\Omega)
 \le
 P_s(F,\Omega)+P_s(C,\Omega).
\end{equation}
Combining \eqref{eq:summary-two-sided} and \eqref{eq:summary-submod} gives
\[
 2P_s(C,\Omega)
 \le
 P_s(F,\Omega)+P_s(C,\Omega),
\]
hence
\[
 P_s(C,\Omega)\le P_s(F,\Omega).
\]
\end{proof}

\end{document}